\documentclass[10pt]{article}
\usepackage{dcolumn}
\usepackage{bm}
\usepackage{verbatim}       

\usepackage[margin=4.3cm]{geometry}

\usepackage{amsthm}

\usepackage{bm}
\usepackage{amstext}
\usepackage{amsmath}
\usepackage{mathtools}
\usepackage{amsfonts}
\usepackage{mathrsfs}
\usepackage{wrapfig}

\newcommand{\p}{\partial}

\newcommand{\bd}{\begin{definition}}                
\newcommand{\ed}{\end{definition}}                  
\newcommand{\bc}{\begin{corollary}}                 
\newcommand{\ec}{\end{corollary}}                   
\newcommand{\bl}{\begin{lemma}}                     
\newcommand{\el}{\end{lemma}}                       
\newcommand{\bp}{\begin{proposition}}            
\newcommand{\ep}{\end{proposition}}                
\newcommand{\bere}{\begin{remark}}                  
\newcommand{\ere}{\end{remark}}                     

\newcommand{\bt}{\begin{theorem}}
\newcommand{\et}{\end{theorem}}

\newcommand{\be}{\begin{equation}}
\newcommand{\ee}{\end{equation}}

\newcommand{\bit}{\begin{itemize}}
\newcommand{\eit}{\end{itemize}}
\newtheorem{theorem}{Theorem}[section]
\newtheorem{corollary}[theorem]{Corollary}
\newtheorem{lemma}[theorem]{Lemma}
\newtheorem{proposition}[theorem]{Proposition}
\theoremstyle{definition}
\newtheorem{definition}[theorem]{Definition}
\theoremstyle{remark}
\newtheorem{remark}[theorem]{Remark}

\begin{document}

\title{Indicatrix geometry clarifies that Finsler length can be larger than  relative length}

\author{E. Minguzzi\thanks{
Dipartimento di Matematica e Informatica ``U. Dini'', Universit\`a
degli Studi di Firenze, Via S. Marta 3,  I-50139 Firenze, Italy.
E-mail: ettore.minguzzi@unifi.it} }


\date{}

\maketitle

\begin{abstract}
\noindent  I show that  Matsumoto conjectured inequality between relative length and Finsler length is false. The incorrectness of the  claim is easily inferred from the geometry of the indicatrix.
\end{abstract}

\section{Introduction}
In Finsler geometry we are given a manifold $M$, and a $C^2$ positive homogeneous function $F\colon TM\backslash 0\to [0,+\infty)$ such that the  vertical Hessian
\[
g_{ab}(x,y)=\frac{1}{2} \frac{\p^2 F^2}{\p y^a \p y^b}
\]
is positive definite, where $(x^a,y^a)$ are the canonical coordinates on $TM$.

The indicatrix  $\mathscr{I}_x$ is the locus $F(x,y)=1$. Since in this work we are just concerned with the geometry of the indicatrix as an embedded submanifold of $V=T_xM$, we could suppress the variable $x$ from all the next formulas.\\

\begin{figure}{h}
\centering
\includegraphics[width=3.5cm]{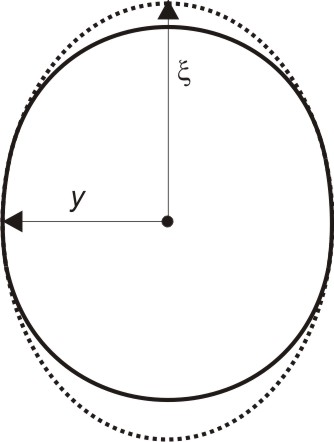}
\caption{The indicatrix $\mathscr{I}_x$ (solid line) and the osculating ellipsoid at $y$ (dotted line).}
\label{nos}
\end{figure}

The Finsler metric $g$ has a very simple geometrical interpretation in terms of the indicatrix \cite{laugwitz11}, see Fig.\ \ref{nos}. Namely for $y\in \mathscr{I}_x$
\[
\mathscr{E}(y)=\{\eta \in V\colon g_{ab}(x,y) \eta^a \eta^b=1 \},
\]
is the osculating ellipsoid to the indicatrix at $y\in \mathscr{I}_x$. The condition of positive definiteness on $g$ is equivalent to the condition of strong convexity for the indicatrix at every point (strong convexity means that, in affine coordinates on the tangent space,  the Taylor expansion of the indicatrix graphing function at every point of the indicatrix has non-trivial second order terms).

Clearly there exists a closed convex set $\mathscr{B}_x$ with strongly convex boundary $\mathscr{I}_x$ whose osculating ellipsoid $\mathscr{E}(y)$ at a point $y$ is not entirely contained in $\mathscr{B}_x$.
Pick a point $\xi\in \mathscr{E}(y)\backslash \mathscr{B}_x$, then $F(x,\xi)>1$ while $g_{ab}(x,y)\xi^a\xi^b=1$, thus defined the {\em relative length} (i.e.\ relative with respect to $y$)
\[\vert \xi\vert_y:=\sqrt{g_{ab}(x,y)\xi^a\xi^b}
\]
we have $\vert \xi\vert_y<F(x,\xi)$.

In \cite{matsumoto79} Matsumoto conjectured that for every $\xi,y\in T_xM\backslash 0$
 the inequality $\vert \xi\vert_y \ge F(x,\xi) $ holds true. Our considerations above show that this is not necessarily the case.
 Of course, a similar argument can be used to show that the inequality $\vert \xi\vert_y \le F(x,\xi) $ for every $\xi,y\in T_xM\backslash 0$, does not necessarily hold, it is sufficient to consider an indicatrix which admits an osculating ellipsoid which does not fully contain the indicatrix.

\end{document}